\documentclass[12pt,reqno]{amsart}

\usepackage{color}

\usepackage[colorlinks=true,
linkcolor=webgreen,
filecolor=webbrown,
citecolor=webgreen]{hyperref}
\definecolor{webgreen}{rgb}{0,.5,0}
\definecolor{webbrown}{rgb}{.6,0,0}

\usepackage{longtable}

\def\trace{\mbox{trace}}

\newtheorem{theorem}{Theorem}[section]
\newtheorem{lemma}[theorem]{Lemma}

\theoremstyle{definition}

\newtheorem{definition}[theorem]{Definition}

\newtheorem{example}[theorem]{Example}

\newtheorem{remark}[theorem]{Remark}

\begin{document}

\title{Arithmetic and growth of periodic orbits}
\author{Yash Puri}
\author{Thomas Ward}
\address{School of Mathematics\\
University of East Anglia\\
Norwich NR4 7TJ\\
U.K.}
\email{\href{mailto:t.ward@uea.ac.uk}{t.ward@uea.ac.uk}}
%\subjclass{}
\date{\today}
\thanks{The first author gratefully acknowledges
the support of E.P.S.R.C. grant 96001638}
\begin{abstract}
Two natural properties of integer
sequences are introduced and studied.
The first, {\sl exact realizability},
is the property that the sequence coincides with the
number of periodic points under some map.
This is shown to impose a strong inner
structure on the sequence. The second,
{\sl realizability in rate}, is the property that
the sequence asympototically approximates the
number of periodic points under some map.
In both cases we discuss when a sequence
can have that property.
For exact realizability, this amounts to
examining the range and domain among integer sequences
of the paired transformations
$$
\mbox{Per}_n=\sum_{d\vert n}d\mbox{Orb}_d;\quad\quad\quad
\mbox{Orb}_d=\frac{1}{n}\sum_{d\vert n}\mu(n/d)\mbox{Per}_d
\quad\quad\quad{\bf ORBIT}
$$
that move between an arbitrary sequence of non-negative
integers Orb counting the orbits of a map and the sequence
Per of periodic points for that map.
Several examples
from the
\htmladdnormallink{Encyclopedia of Integer
Sequences}{http://www.research.att.com/~njas/sequences/}
arise in this work, and a table
of
sequences from the Encyclopedia known or conjectured to be
exactly realizable is given.
\end{abstract}
\maketitle

\tableofcontents

\section{Introduction}

Let $T:X\to X$ be a map.
Three measures of growth in complexity for $T$ are given by
the {\sl number of points with period $n$},
\begin{equation*}
f_n(T)=\#\{x\in X\mid T^nx=x\},
\end{equation*}
the {\sl number of points with least period $n$},
\begin{equation*}
f^{*}_n(T)=\#\{x\in X\mid
T^n(x)=x\mbox{ and }\#\{T^kx\}_{k\in\mathbb Z}=n\},
\end{equation*}
and the {\sl number of orbits of length $n$},
\begin{equation*}
f^{o}_n(T)=f^{*}_n(T)/n.
\end{equation*}
In this note we assume that $f_n(T)$ is finite for $n\ge1$
and give some results on what arithmetic properties the
sequence $(f_n(T))$ may have, and show when the growth
in $(f_n(T))$ is related to the growth in $(f^{*}_n(T))$.
It will be convenient to adopt the following notation:
a sequence $a_1,a_2,a_3,\dots$ is denoted $(a_n)$ or
simply $a$.

\begin{definition}\label{realises}\rm
Let $\phi=(\phi_n)$ be a sequence of non-negative integers. Then
\begin{enumerate}
\item $\phi\in\mathcal{ER}$ ({\sl exactly realizable}) if there is a
set $X$ and a map $T:X\to X$ for which
$f_n(T)=\phi_n$ for all $n\ge 1$;
\item $\phi\in\mathcal{RR}$ ({\sl realizable in rate}) if there is a
set $X$ and a map $T:X\to X$ for which
$f_n(T)/\phi_n\to 1$ as $n\to \infty$.
\end{enumerate}\end{definition}

None of
the results below are changed if the realizing maps
are required to be homeomorphisms of a compact $X$,
but this is not pursued here.

\section{Exact realization}

The set of points with period $n$ under $T$ is the
disjoint union of the set of points with least period
$d$ under $T$ for $d$ dividing $n$, so
\begin{equation}\label{sumoverdivisors}
f_n(T)=\sum_{d\vert n}f^{*}_d(T).
\end{equation}
Equation (\ref{sumoverdivisors}) may be
inverted via the M{\"o}bius inversion formula to give
\begin{equation}\label{basicmobiusinversion}
f^{*}_n(T)=\sum_{d\vert n}\mu(n/d)f_d(T),
\end{equation}
where $\mu(\cdot)$ is the M{\"o}bius function.
On the other hand, the set of points with least
period $n$ comprises exactly $f^{o}_n$ orbits each of
length $n$, so
\begin{equation}\label{divisibilitycondition}
0\le f^{*}_n(T)=\sum_{d\vert n}\mu(n/d)f_d(T)\equiv 0\mbox{ mod $n$}.
\end{equation}
It is clear (since one may take $X=\mathbb N$ and make $T$ to
be a permutation with the appropriate number of cycles of each
length) that these are the only conditions for membership
in $\mathcal{ER}$.

\begin{lemma}\label{exactrealization}
Let $\phi$ be a sequence of non-negative integers.
Then $\phi\in{\mathcal{ER}}$ if and only if
$\sum_{d\vert n}\mu(n/d)\phi_d$ is non-negative and
divisible by $n$ for all $n\ge 1$.
\end{lemma}

Everything that follows is a consequence of this lemma.
Before considering properties of $\mathcal{ER}$ as a whole,
some examples are considered. The sequences that arise here are therefore
close in spirit to the `eigen-sequences' for the
transformation {\bf M{\"O}BIUS} discussed
in \cite{MR96i:05004}
with the additional requirement that the sequence
$f^*$ be divisible by $n$ and non-negative.

\begin{example}\label{firstexamples}\rm
\begin{enumerate}\item The Fibonacci sequence
\htmladdnormallink{A000045}{http://www.research.att.com:80/cgi-bin/access.cgi/as/njas/sequences/eisA.cgi?Anum=000045}
is
not in $\mathcal{ER}$. Using
(\ref{divisibilitycondition}) we see that
$f_3-f_1$ must always be divisible by $3$, but the
Fibonacci sequence begins $1,1,2,3,\dots$.
By contrast the golden mean shift
(see \cite{MR97a:58050}) shows that the
closely related Lucas sequence
\htmladdnormallink{A000204}{http://www.research.att.com:80/cgi-bin/access.cgi/as/njas/sequences/eisA.cgi?Anum=000204}
is in $\mathcal{ER}$.
This will be dealt
with in greater generality in Section \ref{recurrence}
below.

\item For any map $T$,
equation (\ref{divisibilitycondition}), when $n$ is
a prime $p$, states that
$$
f_p(T)\equiv f_1(T)\mbox{ mod $p$}.
$$
If $A\in GL_k(\mathbb Z)$ is an invertible integer
matrix with no unit root eigenvalues, then the
periodic points in the corresponding automorphism of the
$k$-torus show that
$$
\det(A^p-I)\equiv\det(A-I)\mbox{ mod $p$}
$$
for all primes $p$.

\item Similarly, if $B\in M_k(\mathbb N)$ is
a matrix of non-negative integers, the associated subshift of finite
type (see \cite{MR97a:58050}) shows that
$$
\trace(B^p)\equiv\trace(B)\mbox{ mod $p$}
$$
for all primes $p$.
When $k=1$ this is Fermat's little theorem.
When $B=[2]$, so $f_n=2^n$, $f^{o}_n$ is the sequence
\htmladdnormallink{A001037}{http://www.research.att.com:80/cgi-bin/access.cgi/as/njas/sequences/eisA.cgi?Anum=001037}
(shifted by one)
counting
irreducible polynomials of degree $n$ over $\mathbb F_2$.

\item The subshifts of finite type give a
family of elements of $\mathcal{ER}$ of exponential
type. Another family comes from Pascal's triangle:
if $k>1$, $1\le j<k$ and $a_n=\binom{kn}{jn}$, then
$a\in\mathcal{ER}$. For $k=2$ and $j=1$, if $f_n=f_n(T)$ for
the realizing map $T$, then
$f^{*}_n$ is the sequence
\htmladdnormallink{A007727}{http://www.research.att.com:80/cgi-bin/access.cgi/as/njas/sequences/eisA.cgi?Anum=007727}
counting
$2n$-bead black and white strings
with $n$ black beads and fundamental period $2n$.

\item Connected $S$-integer dynamical systems (see
\cite{MR99b:11089}, \cite{MR99k:58152} for these
and the next example):
a subset $S\subset\{2,3,5,7,11,\dots\}$ and a
rational $\xi\neq0$ are given with the property that
$\vert\xi\vert_p>1\implies p\in S$.
The resulting system constructs a map $T:X\to X$
for which
$$
f_n(T)=\prod_{p\le\infty}\vert\xi^n-1\vert_p.
$$
With $\xi=2,S=\{2,3,5,7\}$ this gives the
sequence
$$
1,1,1,1,31,1,127,17,73,341,2047,13,8191,5461,4681,\dots
$$
in $\mathcal{ER}$.

\item Zero-dimensional $S$-integer dynamical
systems: a prime $p$ is fixed, a subset
$S$ of the set of all irreducible polynomials
in $\mathbb F_p[t]$ and a rational function $\xi\in\mathbb F_p(t)$
are given, with the property that
$\vert\xi\vert_{f}>1\implies f\in S$.
The resulting system constructs a map
$T:X\to X$ for which
$$
f_n(T)=\vert\xi^n-1\vert_{t^{-1}}\times\prod_{f\in S}\vert\xi^n-1\vert_f
$$
where $\vert\cdot\vert_{t^{-1}}$ is used to denote
the valuation `at infinity' induced by $\vert t\vert_{t^{-1}}=p$.
Taking $p=2$, $S=\{t-1\}$ and $\xi=t$ gives the
formula
$$
f_n(T)=2^{n-2^{\rm{ord}_2(n)}}
$$
and the
sequence
\htmladdnormallink{A059991}{http://www.research.att.com:80/cgi-bin/access.cgi/as/njas/sequences/eisA.cgi?Anum=059991}
in $\mathcal{ER}$.
\end{enumerate}
\end{example}

\subsection{Algebra of exactly realizable sequences}
The set $\mathcal{ER}$ -- or the ring $K_0(\mathcal{ER})$ --
has a very rich structure. Say that a sequence $a\in\mathcal{ER}$
{\sl factorizes} if there exists sequences $b,c\in\mathcal{ER}$
with $a_n=b_nc_n$ for all $n\ge 1$, and is {\sl prime}
if such a factorization requires one of $b$ or $c$ to be the
constant sequence $(1)$.

\begin{lemma}\label{basicalgebra}
$\mathcal{ER}$ contains the constant sequences
and is closed under addition and multiplication.
Elements of $\mathcal{ER}$ may have infinitely many
non-trivial factors.
There are non-trivial primes in $\mathcal{ER}$.
\end{lemma}

\begin{proof} The constant sequence $(1)$ is in $\mathcal{ER}$
since it is realized by taking $X$ to be a singleton.
The condition in Lemma \ref{exactrealization}
is closed under addition.
On the other hand, if $\phi$ and $\psi$
are exactly realized by systems $(X,T)$
and $(Y,S)$, then $(X\times Y,T\times S)$
exactly realizes $(\phi_n\cdot\psi_n)$.
For each $k\ge 1$ define a sequence
$r^{(k)}$ by $r^{(k)}_n=0$ for
$1<n\le k$ and $r^{(k)}_n=1$ for $n>k$ or $n=1$.
Then $a^{(k)}\in{\mathcal{ER}}$, where
$a_n^{(k)}=\sum_{d\vert n}dr_d^{(k)}$.
Since for each $n$ the sequence
$a_n^{(1)},a_n^{(2)},a_n^{(3)},\dots$ has only finitely
many terms not equal to $1$, the product
$\prod_{k\ge 1}a^{(k)}=(1,3,16,245,1296,41160,\dots)$
is an element of $\mathcal{ER}$ with
infinitely many non-trivial factors.
Finally, the sequence
$(1,3,1,3,\dots)$ is a non-trivial prime in $\mathcal{ER}$.
\end{proof}

In \cite[Sect. 6]{MR86c:58092} a periodic point counting
argument is used to show that the full $p$-shift,
for $p$ a prime, is not topologically conjugate to
the direct product of two dynamical systems. In that argument,
special properties of subshifts of finite type are needed
(specifically, the fact that $f_{p^k}(T)=1$ for all $k\ge 1$
implies that $f_n(T)=1$ for all $n\ge1$ for such systems).
This result does not follow from the arithmetic of
$\mathcal{ER}$ alone: for example, $(3^n)\in{\mathcal{ER}}$ factorizes into
$(1,3,1,3,\dots)\times(3,3,3^3,3^3,\dots)$ in $\mathcal{ER}$ (neither
of which can be realized using a subshift of finite type).
A similar factorization is possible for $(p^n)$ and any odd prime
$p$ (see \cite{puri-thesis} for the details).

\begin{lemma}
There are no non-constant polynomials in $\mathcal{ER}$.
There are non-trivial multiplicative sequences
in $\mathcal{ER}$, but there are
no completely multiplicative sequences apart from
the constant sequence $(1)$.
\end{lemma}

\begin{proof}
Assume that
$$
P(n)=c_0+c_1n+\dots+c_kn^k
$$
with $c_k\neq0$, $k\ge1$,
and that $(P(n))\in{\mathcal{ER}}$.
After multiplying the divisibility condition
(\ref{divisibilitycondition}) by the least
common multiple of the denominators of the
(rational) coefficients of $P$, we
produce a polynomial with integer coefficients
satisfying (\ref{divisibilitycondition}).
It is therefore enough to assume that the
coefficients $c_i$ are all integers.
Let $(f_n)$ and $(f^{*}_n)$ be the periodic
points and least periodic points in the
corresponding system $(X,T)$, and let
$p$ be any prime.
By (\ref{basicmobiusinversion}),
$$
f^{*}_{p^2}=f_{p^2}-f_{p},
$$
so
\begin{eqnarray*}
f^{o}_{p^2}=\frac{f^{*}_{p^2}}{p^2}&=&\frac{f_{p^2}-f_{p}}{p^2}\\
&=&\frac{1}{p^2}
\left(
c_1p^2+c_2p^4+\dots+c_kp^{2k}-
(c_1p+c_2p^2+\dots+c_kp^k)\right)\\
&\in&-\frac{c_1}{p}+\mathbb Z,
\end{eqnarray*}
and therefore $p$ divides $c_1$ for all primes $p$,
showing that $c_1=0$.

Now let $q$ be another prime, and recall that
$$
\mu(1)=1,\mu(p)=-1,\mu(q)=-1,\mu(p^2)=0,
\mu(p^2q)=0,\mu(pq)=1.$$
Since $c_1=0$,
\begin{equation}\label{thatcoolthing}
f_n=c_0+n^2(c_2+c_3n+\dots+c_kn^{k-2}),
\end{equation}
and by (\ref{divisibilitycondition})
$$
p^2q{\big\vert}f_{p^2q}-f_{pq}-f_{p^2}+f_{p}=f^{*}_{p^2q},
$$
so
$$
\frac{f_{p^2}-f_{p}}{p^2q}\in\mathbb Z
$$
by (\ref{thatcoolthing}).
It follows that
$$
c_0(1-1)+c_2(p^4-p^2)+c_3(p^6-p^3)+\dots+
c_k(p^{2k}-p^{k})\in q\mathbb Z
$$
for all primes $q$ and $p$
(since $f_{p^2}-f_{p}$ is certainly divisible by
$p^2$).
So
$$
c_0(1-1)+c_2(p^4-p^2)+c_3(p^6-p^3)+\dots+
c_k(p^{2k}-p^{k})=0;
$$
taking the limit as $p\to\infty$ of
$\frac{1}{p^{2k}}(f_{p^2}-f_{p})$ shows that $c_k=0$.
This contradiction proves the first statement.

There
are many multiplicative sequences in ${\mathcal{ER}}$:
if $f^{*}$ is any multiplicative sequence,
then so is the corresponding
sequence $f$ (see \cite[Theorem 265]{MR81i:10002}).
A multiplicative
sequence $\phi$ is {\sl completely multiplicative}
if
$\phi_{nm}=\phi_{n}\phi_{m}$ for all $n,m\ge 1$.
Assume that $\phi\in{\mathcal{ER}}$ is completely multiplicative,
with $f$ the realising sequence.
For $p$ a prime and any $r\ge 1$,
$$
f^{*}_{p^r}=f_{p^r}-f_{p^{r-1}}=f_p^r-f_p^{r-1}
$$
by (\ref{basicmobiusinversion}).
It follows that
$$
p^r{\big\vert}
f_p^{r-1}(f_p-1).
$$
With $r=1$ this implies that $f_{p}=1+pk_p$ for all $p$, $k_p\in\mathbb N_0$.
Now
$$
p^r{\big\vert}
(1+pk_p)^{r-1}pk_p
$$
for all $p$ and $r\ge 1$. It follows that
$k_p\equiv 0$ mod $p^r$ for all $r\ge 1$, so
$k_p=0$ for all $p$. It follows that
$f_{p}=1$ for all primes $p$, so $f_n=1$ for
all $n\ge 1$.
\end{proof}

Examples show that the
additive convolution
$(\sum_{i+j=n+1,1\le i,j\le n}\phi_i\psi_j)$
of sequences $\phi,\psi\in\mathcal{ER}$
is not
in general in $\mathcal{ER}$.
Similarly, the multiplicative convolution
$(\sum_{d\vert n}\phi_d\psi_{n/d})$ is not in general in
$\mathcal{ER}$.
There is also no closure under quotients:
$(2^n)\in{\mathcal{ER}}$ is term-by-term divisible by
the constant sequence
$(2)\in{\mathcal{ER}}$, but
$(2^{n-1})\notin{\mathcal{ER}}$.

\subsection{Binary recurrence sequences}\label{recurrence}
In this section we expand on the observation made
in Example \ref{firstexamples}.1 by showing that
${\mathcal{ER}}$ only
contains special binary recurrences.

\begin{theorem}\label{binarytheorem}
If $\Delta=a^2+4b$ is not a square, and
$(a,a^2+2b)=1$,
then
a sequence
$u$ with $u_1,u_2\ge1$ satisfying the recurrence
\begin{equation}
\label{basicrecurrence}
u_{n+2}=au_{n+1}+bu_n\mbox{ for }n\ge1
\end{equation}
is in
${\mathcal{ER}}$ if and only
if $\frac{u_2}{u_1}=\frac{a^2+2b}{a}$.
\end{theorem}

As an application, Example \ref{firstexamples}.1
becomes the sharper result that the
Lucasian sequence
$a,b,a+b,a+2b,2a+3b,\dots$ lies in
$\mathcal{ER}$ if and only if $b=3a$.
Moreover, if $f_1=1,f_2=1,f_3=2,\dots$ is
the Fibonacci sequence, then an easy
consequence of Theorem \ref{binarytheorem}
is that for any $k\ge1$ the
sequence $f_k,f_{k+1},f_{k+2},\dots$ is
not in $\mathcal{ER}$.
The more general case with square discriminant,
$a$ and $a^2+2b$ having a common factor
and arbitrary $u_1,u_2$ is dealt with
in \cite{puri-thesis}.

\begin{proof}
First assume that $\frac{u_2}{u_1}=\frac{a^2+2b}{a}$.
Then, by the assumption, the sequence $u$ is
a multiple of the sequence
$a,a^2+2b,a^3+3ab,\dots$ which is
in $\mathcal{ER}$ because the subshift of finite type
corresponding to the matrix
$\bmatrix a&b\cr1&0\endbmatrix$ realizes it (and therefore any multiple of it).

Conversely, assume that $u$ is a sequence in $\mathcal{ER}$
satisfying (\ref{basicrecurrence}). Write $x$ for the
sequence
$$
x:2b,2ab,2(a^2b+b^2),\dots
$$
and $y$ for the
sequence
$$
y:2ab,2(a^2b+2b^2),\dots,
$$
both satisfying the
recurrence (\ref{basicrecurrence}). Notice that
$$
2bu_n=Ax_n+By_n,
$$
for integers $A$ and $B$.
By
(\ref{divisibilitycondition}), for any prime $p$
\begin{equation}\label{cricket}
Ax_p+By_p\equiv Ax_1+By_1\mbox{ mod }p.
\end{equation}
On the other hand, it is well-known that
$x_p\equiv 2b(\frac{\Delta}{p})$ mod $p$ (where
$(\frac{\Delta}{p})$ is the Legendre symbol),
and $y_p\equiv 2ab$ mod $p$ (by the previous
paragraph: $y$ is in $\mathcal{ER}$).
So (\ref{cricket}) implies that
\begin{equation}\label{bat}
2bA\left(\left(\frac{\Delta}{p}\right)-1\right)\equiv 0\mbox{ mod }p
\end{equation}
for all primes $p$.

We now claim that the Legendre symbol $(\frac{\Delta}{p})$ is
$-1$ for infinitely many values of the prime $p$.
This completes the proof of Theorem \ref{binarytheorem},
since (\ref{bat}) forces $A=0$ and hence $u$ is a multiple
of $\frac{1}{2b}y$, namely $a,a^2+2b,\dots$.

To see the claim, choose $c$ such that
$(c,\Delta)=1$ and the Jacobi symbol
$(\frac{c}{\Delta})=-1$. Then by Dirichlet, there are
infinitely many primes $p$ with $p\equiv c$ mod $\Delta$
and $p\equiv 1$ mod $4$.
For such primes,
$(\frac{p}{\Delta})=(\frac{\Delta}{p})=-1$,
which completes the proof.
\end{proof}

The case of square discriminant is much more
involved. A full treatment is in
\cite{puri-thesis}; here we simply show by
examples that the result as stated no longer holds in
general.

\begin{example}\rm
\label{squarediscexample}
\begin{enumerate}
\item There are infinitely many
possible values of the ratio $\frac{u_2}{u_1}$
for binary recurrent sequences in ${\mathcal{ER}}$
satisfying
\begin{equation}\label{squarediscrecurrenceexample}
u_{n+2}=u_{n+1}+2u_n.
\end{equation}
To see this we construct two different realizing
examples and then take linear integral combinations
of them.
The first is the subshift of finite type
$T$ corresponding to the matrix
$A=\bmatrix 1&2\\1&0\endbmatrix$.
This system has (by \cite[Proposition 2.2.12]{MR97a:58050})
$f_n(T)=\trace(A^n)$, which is the
sequence of Jacobsthal--Lucas numbers
\htmladdnormallink{A014551}{http://www.research.att.com:80/cgi-bin/access.cgi/as/njas/sequences/eisA.cgi?Anum=014551}:$1,
5,7,17,\dots$ (shifted by one) and has
initial ratio ${5}$.
On the other hand, the algebraic dynamical system $S$
dual to $x\mapsto -2x$ on the discrete group
$\mathbb Z[\frac{1}{2}]$ has (see, for example,
\cite[Lemma 5.2]{MR99b:11089})
$f_n(S)=\vert(-2)^n-1\vert$, which begins
$3,3,9,15,\dots$ and has ratio $1$.
Now we may apply Lemma \ref{basicalgebra}
as follows. If $s,t\in\mathbb N$ then
$\left(tf_n(T)+sf_n(S)\right)$
in ${\mathcal{ER}}$
is a sequence satisfying (\ref{squarediscrecurrenceexample}).
It follows that the set of possible ratios
$\frac{u_2}{u_1}$ contains the infinite set
$\{\frac{5t+3s}{t+3s}\mid s,t\in\mathbb N\}$.
\item A simpler example
is given by
the Mersenne recurrence.
Since $(2^n)$ and $(1)$ are both in $\mathcal{ER}$,
for any $t,s\ge0$ the sequence $\left(t2^n+s\right)$
satisfying the recurrence
\begin{equation}\label{mersenne}
u_{n+2}=3u_{n+1}-2u_n
\end{equation}
is in $\mathcal{ER}$. Thus the
set of possible ratios $\frac{u_2}{u_1}$
for exactly realizable solutions of (\ref{mersenne})
contains the infinite set
$\left\{\frac{4t+s}{2t+s}\mid s,t\in\mathbb N\right\}.$
\end{enumerate}
\end{example}

For higher order recurrences with
companion polynomials irreducible over
the rationals, it is clear that some
analogue of Theorem \ref{binarytheorem} holds. The rational
solutions of a $k$th order reccurence form a rational
$k$-space; the smallest subspace contained in $\mathcal{ER}$
has dimension strictly smaller than $k$. Is this dimension
always 1?

\section{Realization in rate}

Write $\lfloor x\rfloor$ for the greatest integer
less than or equal to $x$ and $\lceil x\rceil$ for the
smallest integer greater than or equal to $x$.
In this section we assume that sequences are
never zero. Different complications arise from
zeros of sequences and these are discussed in
detail in \cite{puri-thesis}.

\begin{theorem}\label{basicrateresults} Let $\alpha,\beta$ be
positive constants.\begin{enumerate}
\item If $\phi_n\to\infty$ with $\frac{\phi_n}{n}\to 0$, then
$\phi\notin{\mathcal{RR}}$.
\item The sequence $(\lfloor n^{\alpha}\rfloor)\in\mathcal{RR}$ if
and only if $\alpha>1$.
\item The sequence $(\lfloor\beta^n\rfloor)\in\mathcal{RR}$ if and
only if $\beta\ge1$.
\end{enumerate}
\end{theorem}

\begin{proof}
\noindent1. Assume that $\phi\in{\mathcal{RR}}$ and let $f$ be the
corresponding sequence of periodic points.
Then $\frac{f_n}{\phi_n}\to 1$, so
$\{\frac{f_n}{\phi_n}\}$ is bounded.
It follows that
$\{\frac{f^{*}_n}{\phi_n}=\frac{n}{\phi_n}f^{o}_n\}$ is bounded, and
hence
$f^{*}_n=nf^{o}_n=0$ for all large $n$. This implies that $f_n$ is bounded,
and so $\frac{f_n}{\phi_n}\to0$, which
contradicts the assumption.

\noindent2. For $\alpha\in(0,1)$ this follows from
part 1.
Suppose therefore that $(n)\in{\mathcal{RR}}$.
Then there is a sequence $f\in\mathcal{ER}$
with $f_n/n\to 1$, so for $p$ a prime,
$pf^{o}_p=f^{*}_p=f_p-f^{*}_1$, and therefore
$f^{o}_p\to 1$ as $p\to\infty$.
Since $f^{o}_p$ is an integer, it follows that $f^{o}_p=1$ for
all large $p$.
Now let $q$ be another large prime.
Then
$$
\frac{f_{pq}}{pq}=\frac{f^{*}_{pq}+f^{*}_p+f^{*}_q+f^{*}_1}{pq}=
\frac{f^{*}_{pq}}{pq}+\frac{1}{p}+\frac{1}{q}+\frac{f^{*}_1}{pq},
$$
so
$$
\frac{1}{p}+\frac{1}{q}+\frac{f^{*}_1}{pq}-\frac{f_{pq}}{pq}\in\mathbb Z.
$$
Fix $p$ large and let $q$ tend to infinity to see that
$$
\frac{1}{p}\in\mathbb Z,
$$
which is impossible.
The same argument shows that $f_n/n$ cannot
have any positive limit as $n\to\infty$.

For $\alpha>1$, let $f^{o}_n=\lceil n^{\alpha-1}\prod_{p\vert d}
(1-p^{-\alpha})\rceil$, where the product runs over prime
divisors only.
Then
$$
\sum_{d\vert n}d^{\alpha}\prod_{p\vert d}(1-p^{-\alpha})=n^{\alpha}
\le
\sum_{d\vert n}df^{o}_d=f_n\le n^{\alpha}+\sum_{d\vert n}d,
$$
so $0\le f_n-\phi_n\le o(n^{\alpha}).$

\noindent3. This is clear: for $\beta<1$ the sequence is
eventually $0$; for $\beta>1$ the construction used
in part 2. works.
\end{proof}

There are sequences growing more slowly than $n^{\alpha}$ in
$\mathcal{RR}$: in \cite[Chap. 5]{puri-thesis} it is shown that
$\left(\lfloor Cn^s(\log n)^r\rfloor\right)\in\mathcal{RR}$
for any $r\ge 1,C>0,s\ge1$.

\section{Comparing orbits with periodic points}

As is well-known, if $f^{*}$ grows fast enough,
then $f$ grows very much like $f^{*}$ (though
not conversely in the case of super-exponential growth:
cf. Theorem \ref{superexponential} below).
Throughout this section $f_n=f_n(T)$ and $f^{*}_n=f^{*}_n(T)$ for some
map $T$.

\begin{remark}\rm That $f^{*}_n$ is close to $f_n$ when $f_n$ is growing
exponentially has been commented on by
Lind in \cite[Sect. 4]{MR84g:28017}.
He points out, using
(\ref{basicmobiusinversion}), that if $T$ is the
automorphism of the $2$-torus corresponding to the
matrix
$\bmatrix 2&1\\1&1\endbmatrix$
then
$f^{*}_{20}(T)$ is only $0.006\%$ smaller than
$f_{20}(T)$. The sequence $f$ of periodic points for
this map is
\htmladdnormallink{A004146}{http://www.research.att.com:80/cgi-bin/access.cgi/as/njas/sequences/eisA.cgi?Anum=004146}.
\end{remark}

\begin{theorem}\begin{enumerate}\label{superexponential}
\item If $\frac{1}{n}\log f^{*}_n\to C\in[0,\infty]$ then
$\frac{1}{n}\log f_n\to C$ also.
\item $\frac{1}{n}\log f^{*}_n\to C\in(0,\infty)$ if and
only if $\frac{1}{n}\log f_n\to C$.
\item If $\frac{1}{n}\log f_n\to\infty$ then
$\left\{\frac{1}{n}\log f^{*}_n\right\}$ may be unbounded
with infinitely many limit points.
\end{enumerate}
\end{theorem}

\begin{proof}
1. If $\frac{1}{n}\log f^{*}_n\to\infty$ then
$\frac{1}{n}\log f_n\to\infty$ also, since
$f_n\ge f^{*}_n$ for all $n$.
If $\frac{1}{n}\log f^{*}_n\to C\in[0,\infty)$, then
(for $n$ large enough to have $f^{*}_n\neq 0$)
\begin{eqnarray*}
\frac{1}{n}\log f^{*}_n\le\frac{1}{n}\log f_n&=&
\frac{1}{n}\log\left(
\sum_{d\vert n}f^{*}_d\right)\\
&\le&\frac{1}{n}\log n+\frac{1}{n}\log\max_{d\vert n}
\{f^{*}_d\}.
\end{eqnarray*}
For each such
$n$, choose $\tilde{n}\in\{d\mid
d\vert n, f^{*}_d\ge f^{*}_{d'}{\ }\forall{\ }d'\vert n\}$
so that $f^{*}_{\tilde{n}}=\max_{d\vert n}\{f^{*}_d\}$
and $\frac{\tilde{n}}{n}\le 1$.
Then
\begin{eqnarray*}
\frac{1}{n}\log f_n&\le&
\frac{1}{n}\log n+\frac{\tilde{n}}{n}\cdot
\frac{1}{\tilde{n}}\log f^{*}_{\tilde{n}}\\
&\le&
\frac{1}{n}\log n+\frac{1}{\tilde{n}}\log f^{*}_{\tilde{n}}
\to C.
\end{eqnarray*}

\noindent 2. It is enough to show that if
$\frac{1}{n}\log f_n\to C\in(0,\infty)$ then
$\frac{1}{n}\log f^{*}_n\to C$ also.
For $r\ge 1$,
\begin{equation*}
f_r\ge f^{*}_r=-\sum_{d\vert r,d\neq r}f^{*}_d+f_r\ge
f_r-\sum_{d\vert r,d\neq r}f_d.
\end{equation*}
Let $R$ be an upper bound for
$\{\frac{1}{n}\log f_n\mid
f_n\neq0\}$ and pick $\epsilon\in(0,3C)$.
Choose $N$ so that
$$
r>N\implies
e^{r(C-\epsilon)}\le f_r\le
e^{r(C+\epsilon)}.
$$
Then for $r>2N$ (so that
$r^{\ast},N\le\lfloor\frac{r}{2}\rfloor$),
\begin{eqnarray*}
f_r\ge f^{*}_r&\ge&f_r-\sum_{n=1}^{N}f_n-\sum_{n=N+1}^{\lfloor r/2\rfloor}f_n\\
&\ge&
f_r-\left(Ne^{NR}+(r/2-N)e^{r(C+\epsilon)/2}\right)\\
&\ge&
f_r\left(1-Ne^{NR-r(C-\epsilon)}-(r/2-N)e^{-r(C-3\epsilon)/2}\right),
\end{eqnarray*}
and the bracketed expression converges to $1$ as $r\to\infty$.
Taking logs and dividing by $r$ gives the result.

\noindent 3. Write $p_1,p_2,\dots$ for the sequence of primes.
Let $n_r=p_rp_{r+1}$, and define a sequence $(f^{*}_k)$
as follows.
For $k$ not of the form $n_r$, define $f^{*}_k=k\cdot 2^{k^3}.$
For $k$ of the form $n_r$ define $f^{*}_k$ according to the
following scheme:
\begin{eqnarray*}
f^{*}_{n_1}&=&n_12^{n_1}\\
f^{*}_{n_2}&=&n_22^{n_2},f^{*}_{n_3}=n_32^{2n_3}\\
f^{*}_{n_4}&=&n_42^{n_4},f^{*}_{n_5}=n_52^{2n_5},f^{*}_{n_6}=n_62^{3n_6}\\
f^{*}_{n_7}&=&n_72^{n_7},f^{*}_{n_8}=n_82^{2n_8},f^{*}_{n_9}=n_92^{3n_9},
f^{*}_{n_{10}}=n_{10}2^{4n_{10}}\\
\end{eqnarray*}
and so on.
Then
$\frac{1}{n}\log f_n\to\infty$ off the
$n_r$'s clearly. Along the sequence $(n_r)$,
$$
f_{n_r}=f^{*}_{n_r}+f^{*}_{p_r}+f^{*}_{p_{r+1}}+f^{*}_1\ge
f^{*}_{p_{r+1}},
$$
so
$$
\frac{1}{n_r}\log f_{n_r}\ge
\frac{1}{p_rp_{r+1}}\log\left(
p_{r+1}\cdot
2^{p_{r+1}^3}\right)\to\infty.
$$

On the other hand, along a subsequence of $n_r$'s
chosen to have
$f^{*}_{n_r}=n_r2^{\ell n_r}$ for a fixed
$\ell\in\mathbb N$ (which will exist
by construction), we realize
$\ell\log 2$ as a limit point of the
sequence $\frac{1}{n}\log f^{*}_n$.
\end{proof}

Finally, we turn to comparing these growth rates
in a sub-exponential setting. For polynomial
growth, the next result shows that $f$ and
$f^{*}$ are forced to behave very differently.

\begin{theorem}\label{polyratescompare}
Let $C$ and $\alpha$ be positive
constants.
\begin{enumerate}
\item For $\alpha>1$,
the set $\{\frac{f^{*}_n}{n^{\alpha}}\}$ is bounded
if and only if $\{\frac{f_n}{n^{\alpha}}\}$ is bounded.
\item For $\alpha>1$, $\frac{f_n}{n^{\alpha}}\to0$
if and only if
$\frac{f^{*}_n}{n^{\alpha}}\to0$.
\item If $\frac{f_n}{n^{\alpha}}\to C$ for some $\alpha>1$,
then $\{\frac{f^{*}_n}{n^{\alpha}}\}$ has infinitely
many limit points.
\item If $\frac{f^{*}_n}{n^{\alpha}}\to C$ for some $\alpha\ge1$,
then $\{\frac{f_n}{n^{\alpha}}\}$ has infinitely
many limit points.
\end{enumerate}
\end{theorem}

\begin{proof} 1. Let $R$ be an upper bound for
$\{\frac{f^{*}_n}{n^{\alpha}}\}$. Then
\begin{equation*}
\frac{f_n}{n^{\alpha}}\le
\frac{1}{n^{\alpha}}\sum_{d\vert n}Rd^{\alpha}
=
R\sum_{d\vert n}\left(\frac{d}{n}\right)^{\alpha}
\le R\sum_{d=1}^{\infty}\frac{1}{d^{\alpha}}<\infty.
\end{equation*}
The converse is obvious.

\noindent2. One direction is clear.
Assume that $\frac{f^{*}_n}{n^{\alpha}}\to 0$.
Fix $\epsilon>0$; choose
$M_1\in\mathbb N$ so that
$$
n>M_1\implies
\frac{f^{*}_n}{n^{\alpha}}<\frac{\epsilon}{1+\beta}
$$
where $\beta=\sum_{k=1}^{\infty}\frac{1}{k^{\alpha}}$.
Choose $M_2$ so that
$$
n>M_2\implies
\sum_{r=1}^{M_1}\frac{f^{*}_r}{n^{\alpha}}<\frac{\epsilon}{1+\beta}.
$$
Then for $n\ge\max\{M_1,M_2\}$,
\begin{eqnarray*}
0\le\frac{f_n}{n^{\alpha}}=\sum_{d\vert n}\frac{f^{*}_d}{n^{\alpha}}
&\le&\sum_{r=1}^{M_1}\frac{f^{*}_r}{n^{\alpha}}+
\sum_{d\vert n,d>M_1}\frac{f^{*}_r}{n^{\alpha}}\\
&=&\sum_{r=1}^{M_1}\frac{f^{*}_r}{n^{\alpha}}+
\sum_{d\vert n,d>M_1}+
\frac{d^{\alpha}}{n^{\alpha}}\cdot\frac{f^{*}_d}{d^{\alpha}}\\
&\le&
\frac{\epsilon}{1+\beta}+\frac{\epsilon}{1+\beta}
\sum_{d\vert n,d>M_1}\frac{d^{\alpha}}{n^{\alpha}}\\
&\le&\frac{\epsilon}{1+\beta}+\beta
\frac{\epsilon}{1+\beta}\le\epsilon.
\end{eqnarray*}

\noindent3. Assume that $\frac{f_n}{n^{\alpha}}\to C>0$.
Then $\frac{f^{*}_p}{p^{\alpha}}\to C$ along primes.
For a fixed prime $p$,
$$
\frac{f^{*}_{p^r}}{p^{r\alpha}}=\frac{f_{p^r}}{p^{r\alpha}}
-\frac{f_{p^{r-1}}}{p^{(r-1)\alpha}}\cdot\frac{1}{p^{\alpha}}
\to\left(1-\frac{1}{p^{\alpha}}\right)C
$$
as $r\to\infty$.

\noindent4. Assume that $\frac{f^{*}_n}{n^{\alpha}}\to C>0$.
Then $\frac{f_p}{p^{\alpha}}\to C$ along primes.
For fixed prime $p$ and $q$ prime,
$$
\frac{f_{pq}}{(pq)^{\alpha}}=\frac{f^{*}_{pq}+f^{*}_{q}+f^{*}_{p}+f^{*}_{1}}{(pq)^{\alpha}}
\to\left(1+\frac{1}{p^{\alpha}}\right)C
$$
as $q\to\infty$.
\end{proof}

\begin{remark}\rm For the case
$\frac{f^{*}_n}{n}\to C>0$ in Theorem \ref{polyratescompare},
$\frac{f_n}{n}$ is unbounded: similar arguments show that
$$
\frac{f_{p_1p_2\dots p_m}}{p_1p_2\dots p_m}\ge
\sum_{i=1}^{m}\frac{1}{p_i}\to \infty
$$
as $m\to\infty$.
\end{remark}

\section{Examples}

Few of the standard sequences turn out to be
in $\mathcal ER$. Here we list a few that are,
and one that nearly is.
In some cases the
proof proceeds by exhibiting a realizing map,
in others by proving the congruence.
Section \ref{summary} contains a table
with many sequences from the Encyclopedia in
$\mathcal{ER}$; in particular all
sequences realized by oligomorphic permutation
groups from \cite{MR1750744} that fall in $\mathcal{ER}$ are
listed.

\begin{example}\begin{enumerate}\rm
\item Many trivial sequences are in $\mathcal{ER}$,
among them
\htmladdnormallink{A00004}{http://www.research.att.com:80/cgi-bin/access.cgi/as/njas/sequences/eisA.cgi?Anum=00004},
\htmladdnormallink{A00012}{http://www.research.att.com:80/cgi-bin/access.cgi/as/njas/sequences/eisA.cgi?Anum=00012},
\htmladdnormallink{A00079}{http://www.research.att.com:80/cgi-bin/access.cgi/as/njas/sequences/eisA.cgi?Anum=00079}
(shifted by one),
\htmladdnormallink{A00203}{http://www.research.att.com:80/cgi-bin/access.cgi/as/njas/sequences/eisA.cgi?Anum=000203">A000203}.
\item \htmladdnormallink{A023890}{http://www.research.att.com:80/cgi-bin/access.cgi/as/njas/sequences/eisA.cgi?Anum=023890},
the sum of non-prime divisors, is
in $\mathcal{ER}$ since it corresponds to having one
orbit of each composite length.
\item \htmladdnormallink{A000984}{http://www.research.att.com:80/cgi-bin/access.cgi/as/njas/sequences/eisA.cgi?Anum=000984}
(shifted
by one). As pointed out in Example \ref{firstexamples}.4,
the sequence of central binomial coefficients
$\binom{2n}{n}$ is in $\mathcal{ER}$ for a combinatorial
reason.
Similarly the sequences of the form
$\binom{kn}{jn}$ are all in $\mathcal{ER}$: these
include
\htmladdnormallink{A005809}{http://www.research.att.com:80/cgi-bin/access.cgi/as/njas/sequences/eisA.cgi?Anum=005809}
($k=3,j=1$).
\item \htmladdnormallink{A001035}{http://www.research.att.com:80/cgi-bin/access.cgi/as/njas/sequences/eisA.cgi?Anum=001035}
(shifted by one)
counts the number of distinct posets on $n$ labeled
elements. The first 16 terms of this sequence
are known, and so the congruence (\ref{divisibilitycondition})
can be verified for $n\le16$. However, the sequence
is not in $\mathcal{ER}$. We are grateful to
Greg Kuperberg for suggesting the following
explanation.
Write $\mathcal P(n)$ for the set of poset
structures on $\mathbb Z/n\mathbb Z$. Then
for $d\vert n$, there is an injection
$\phi_{d,n}:\mathcal P(d)\to\mathcal P(n)$
obtained by pulling back a poset structure
using the 
canonical
homomorphism $\mathbb Z/n\mathbb Z\to\mathbb Z/d\mathbb Z$.
For certain values of $n$, including all prime
values, we claim that those posets that do not appear in
the image of one of these injections
come
in families of size a multiple of $n$, which
gives the congruence (\ref{divisibilitycondition}).
by M{\"o}bius inversion.
Translation gives an action of
$\mathbb Z/n\mathbb Z$ on $\mathcal P(n)$;
if a given poset lies on a free orbit then
that orbit is the family. In general,
suppose that the stabilizer of an orbit
is $\mathbb Z/(n/d)\mathbb Z$, but it is not
in the image of $\phi_{d,n}$.
Then there is a natural action of
the wreath product $\mathbb Z/d\mathbb Z{\rm{wr}}S_{n/d}$
defined by permuting the points in each coset
of $\mathbb Z/(n/d)\mathbb Z$ and adding
a multiple of $n/d$. If $n$ is the
product of two primes (and for many other
$n$) then the size of the orbits
of this action are divisible by $n$.
However, at $n=18$ there are orbits of size
$\pm6$ mod $18$, so here we expect the
congruence (\ref{divisibilitycondition}) to fail.
\item \htmladdnormallink{A001945}{http://www.research.att.com:80/cgi-bin/access.cgi/as/njas/sequences/eisA.cgi?Anum=001945}: $1,
1,1,5,1,7,8,5,19,\dots$
is in $\mathcal ER$ since it counts the
periodic points in the automorphism of the
$3$-torus given by the matrix
$\bmatrix 0&1&0\\0&0&1\\1&1&0\endbmatrix$.
This sequence has been studied computationally
for prime appearances (see \cite{MR1783409}) and
it comes from the cubic polynomial with
smallest Mahler measure (see \cite{MR2000e:11087}).
\item The large class of elliptic divisibility
sequences (see \cite{elliptic}) and Somos sequences seem never to
fall in $\mathcal{ER}$.
\item Three interesting sequences that seem to be
in $\mathcal{ER}$ are the Euler sequence
\htmladdnormallink{A000364}{http://www.research.att.com:80/cgi-bin/access.cgi/as/njas/sequences/eisA.cgi?Anum=000364}
and the sequences
\htmladdnormallink{A006953}{http://www.research.att.com:80/cgi-bin/access.cgi/as/njas/sequences/eisA.cgi?Anum=006953},
\htmladdnormallink{A006954}{http://www.research.att.com:80/cgi-bin/access.cgi/as/njas/sequences/eisA.cgi?Anum=006954}
connected with
the Bernoulli numbers.
\end{enumerate}
\end{example}

\begin{example}Sequences in
$\mathcal{ER}$ arise
from the combinatorics of an iterated map.
It is a natural question to ask what an orbit of
{\bf ORBIT} looks
like, and whether there are any asymptotic properties associated to it.
The simplest orbit starts with the unit sequence
\htmladdnormallink{A000007}{http://www.research.att.com:80/cgi-bin/access.cgi/as/njas/sequences/eisA.cgi?Anum=000007}.
Applying {\bf ORBIT} iteratively
to this gives the following sequence of sequences
(in each case, the sequence counts
the number of periodic points in a map which has the
number of orbits of length $n$ given by the $n$th entry in the previous
sequence).\smallskip

\noindent$
1,0,0,0,0,0,0,0,0,0,0,0,0,0,0,0,0,0,0,0,0,0,0,0,0,0,
\mbox{\hfill(\htmladdnormallink{A000007}{http://www.research.att.com:80/cgi-bin/access.cgi/as/njas/sequences/eisA.cgi?Anum=000007})}
$
\smallskip

\noindent$
1, 1, 1, 1, 1, 1, 1, 1, 1, 1, 1, 1, 1, 1,1,1,1,1,1,1,1,1,1,1,1,1
\mbox{\hfill(\htmladdnormallink{A000012}{http://www.research.att.com:80/cgi-bin/access.cgi/as/njas/sequences/eisA.cgi?Anum=000012})}
$
\smallskip

\noindent$
1, 3, 4, 7, 6, 12, 8, 15, 13, 18, 12, 28, 14, 24, 24, 31,18,39,20
\mbox{\hfill(\htmladdnormallink{A000203}{http://www.research.att.com:80/cgi-bin/access.cgi/as/njas/sequences/eisA.cgi?Anum=000203})}
$
\smallskip

\noindent$
1, 7, 13, 35, 31, 91, 57, 155, 130, 217, 133, 455, 183,399,403
\mbox{\hfill(\htmladdnormallink{A001001}{http://www.research.att.com:80/cgi-bin/access.cgi/as/njas/sequences/eisA.cgi?Anum=001001})}
$
\smallskip

\noindent$
1, 15, 40, 155, 156, 600, 400, 1395, 1210, 2340, 1464, 6200, 2380, 6000 
$
\smallskip

\noindent$
1, 31, 121, 651, 781, 3751, 2801, 11811, 11011, 24211, 16105, 78771, 30941 
$
\smallskip

\noindent$
1, 63, 364, 2667, 3906, 22932, 19608, 97155, 99463, 246078, 177156
$
\smallskip

\noindent$
1, 127, 1093, 10795, 19531, 138811, 137257, 788035, 896260, 2480437 
$
\smallskip

\noindent$
1, 255, 3280, 43435, 97656, 836400, 960800, 6347715, 8069620, 24902280 
$
\smallskip

\noindent$
1, 511, 9841, 174251, 488281, 5028751, 6725601, 50955971, 72636421
$
\smallskip

\noindent$
1, 1023, 29524, 698027, 2441406, 30203052, 47079208, 408345795
$
\smallskip

\noindent$
1, 2047, 88573, 2794155, 12207031, 181308931, 329554457, 3269560515 
$
\smallskip

\noindent$
1, 4095, 265720, 11180715, 61035156, 1088123400, 2306881200 
$
\smallskip

\noindent$
1, 8191, 797161, 44731051, 305175781, 6529545751, 16148168401 
$
\smallskip

\noindent$
1, 16383, 2391484, 178940587, 1525878906, 39179682372, 113037178808 
$
\smallskip

\noindent$
1, 32767, 7174453, 715795115, 7629394531, 235085301451
$
\smallskip

\noindent$
1, 65535, 21523360, 2863245995, 38146972656, 1410533397600 
$
\smallskip

The arithmetic and growth properties of these sequences will be
explored elsewhere. The sequence of first, second and third
terms comprise
\htmladdnormallink{A000012}{http://www.research.att.com:80/cgi-bin/access.cgi/as/njas/sequences/eisA.cgi?Anum=000012},
\htmladdnormallink{A000225}{http://www.research.att.com:80/cgi-bin/access.cgi/as/njas/sequences/eisA.cgi?Anum=000225},
and
\htmladdnormallink{A003462}{http://www.research.att.com:80/cgi-bin/access.cgi/as/njas/sequences/eisA.cgi?Anum=003462} respectively.
\end{example}

\section{Summary}\label{summary}

Being exactly realizable is a strong symmetry property of an
integer sequence. In this table we summarize the
sequences
from the Encyclopedia found to be exactly realizable,
together with the corresponding sequence counting the
orbits, and any
other information. All the sequences are expected to
have realizing maps -- the inclusion of a map means that
we know of a
map that is natural in some sense (for example, has a finite
description or is algebraic). Direct proofs of the congruence are
cited in some brief fashion -- a question mark indicates
that we do not know a proof and seek one, e means it is easy,
and a combinatorial counting problem suggests
why the number of orbits is a non-negative integer.
The combinatorial counting problems and maps are labelled as
follows.\begin{itemize}
\item POLY: the orbits count the number of irreducible polynomials
over a finite field.
\item NECK($k$): the orbits count the number of aperiodic
necklaces with $n$ beads in $k$ colours.
\item NECK: the orbits count a family of necklaces with
constraint -- see the encyclopedia entry for details.
\item KUMMER: follows from the Kummer and von Staudt
congruences.
\item COMB: follows from standard combinatorics arguments.
\item CHK: the orbit sequence is a `CHK' transform.
\item S(1):{$S$-integer map with $\xi=2, S=\{2,3\}, k={\mathbb Q}$}.
\item S(2):$S$-integer map with $\xi=t, S=\{t+1\},
k ={\mathbb F}_2(t)$.
\item R: irrational circle rotation.
\end{itemize}
Of course there are often many ways to fill in the
last column.
If there
is a natural realizing map, then that fact in itself is
usually the best proof of the congruence. Sequences marked with
a
question mark in the first column are not known to be in
$\mathcal{ER}$ at all: they just seem to satisfy the congruence
for the first
twenty or so terms. A star indicates that the initial term
of the sequence is shifted by one. Of course any non-negative
integer sequence at all can appear in the second column,
so the selection here is based on the following arbitrary criterion:
either the periodic point sequence or the orbit counting
sequence is `interesting'. 

\begin{longtable}{llll}\hline
{$f_n(T)$}&{$f_n^o(T)$}&{$T$}
&{Proof of (\ref{divisibilitycondition})}\\
\hline\\
{\htmladdnormallink{A00004}{http://www.research.att.com:80/cgi-bin/access.cgi/as/njas/sequences/eisA.cgi?Anum=00004}}
&
{\htmladdnormallink{A00004}{http://www.research.att.com:80/cgi-bin/access.cgi/as/njas/sequences/eisA.cgi?Anum=00004}}
&{R}&{e}\\
{\htmladdnormallink{A00012}{http://www.research.att.com:80/cgi-bin/access.cgi/as/njas/sequences/eisA.cgi?Anum=00012}}
&
{\htmladdnormallink{A00007}{http://www.research.att.com:80/cgi-bin/access.cgi/as/njas/sequences/eisA.cgi?Anum=00007}}
&
{singleton}&{e}\\
{\htmladdnormallink{A000079*}{http://www.research.att.com:80/cgi-bin/access.cgi/as/njas/sequences/eisA.cgi?Anum=000079}}
&
{\htmladdnormallink{A001037}{http://www.research.att.com:80/cgi-bin/access.cgi/as/njas/sequences/eisA.cgi?Anum=001037}}
&{full 2-shift}&{POLY}\\
{\htmladdnormallink{A000203}{http://www.research.att.com:80/cgi-bin/access.cgi/as/njas/sequences/eisA.cgi?Anum=000203}}
&
{\htmladdnormallink{A00012}{http://www.research.att.com:80/cgi-bin/access.cgi/as/njas/sequences/eisA.cgi?Anum=00012}}
&-&e\\
{\htmladdnormallink{A000204}{http://www.research.att.com:80/cgi-bin/access.cgi/as/njas/sequences/eisA.cgi?Anum=000204}}
&
{\htmladdnormallink{A006206}{http://www.research.att.com:80/cgi-bin/access.cgi/as/njas/sequences/eisA.cgi?Anum=006206}}
&golden mean shift&NECK\\
{\htmladdnormallink{A000244*}{http://www.research.att.com:80/cgi-bin/access.cgi/as/njas/sequences/eisA.cgi?Anum=000244}}
&
{\htmladdnormallink{A027376*}{http://www.research.att.com:80/cgi-bin/access.cgi/as/njas/sequences/eisA.cgi?Anum=027376}}
&{full 3-shift}&{POLY}\\
{\htmladdnormallink{A000302*}{http://www.research.att.com:80/cgi-bin/access.cgi/as/njas/sequences/eisA.cgi?Anum=000302}}
&
{\htmladdnormallink{A027377*}{http://www.research.att.com:80/cgi-bin/access.cgi/as/njas/sequences/eisA.cgi?Anum=027377}}
&{full 4-shift}&{POLY}\\
{\htmladdnormallink{A000351*}{http://www.research.att.com:80/cgi-bin/access.cgi/as/njas/sequences/eisA.cgi?Anum=000351}}
&
{\htmladdnormallink{A001692*}{http://www.research.att.com:80/cgi-bin/access.cgi/as/njas/sequences/eisA.cgi?Anum=001692}}
&{full 5-shift}&{POLY}\\
{\htmladdnormallink{A000364*?}{http://www.research.att.com:80/cgi-bin/access.cgi/as/njas/sequences/eisA.cgi?Anum=000364}}
&
{\htmladdnormallink{A060164}{http://www.research.att.com:80/cgi-bin/access.cgi/as/njas/sequences/eisA.cgi?Anum=060164}}
&{-}&{-}\\
{\htmladdnormallink{A000400*}{http://www.research.att.com:80/cgi-bin/access.cgi/as/njas/sequences/eisA.cgi?Anum=000400}}
&
{\htmladdnormallink{A032164}{http://www.research.att.com:80/cgi-bin/access.cgi/as/njas/sequences/eisA.cgi?Anum=032164}}
&{full 6-shift}&{NECK(6)}\\
{\htmladdnormallink{A000420*}{http://www.research.att.com:80/cgi-bin/access.cgi/as/njas/sequences/eisA.cgi?Anum=000420}}
&
{\htmladdnormallink{A001693}{http://www.research.att.com:80/cgi-bin/access.cgi/as/njas/sequences/eisA.cgi?Anum=001693}}
&{full 7-shift}&{POLY}\\
{\htmladdnormallink{A000593}{http://www.research.att.com:80/cgi-bin/access.cgi/as/njas/sequences/eisA.cgi?Anum=000593}}
&
{\htmladdnormallink{A000035*}{http://www.research.att.com:80/cgi-bin/access.cgi/as/njas/sequences/eisA.cgi?Anum=000035}}
&{-}&{e}\\
{\htmladdnormallink{A000670*}{http://www.research.att.com:80/cgi-bin/access.cgi/as/njas/sequences/eisA.cgi?Anum=000670}}
&
{\htmladdnormallink{A060223}{http://www.research.att.com:80/cgi-bin/access.cgi/as/njas/sequences/eisA.cgi?Anum=060223}}
&{e}\\
{\htmladdnormallink{A000984*}{http://www.research.att.com:80/cgi-bin/access.cgi/as/njas/sequences/eisA.cgi?Anum=000984}}
&
{\htmladdnormallink{A060165}{http://www.research.att.com:80/cgi-bin/access.cgi/as/njas/sequences/eisA.cgi?Anum=060165}}
&{-}&{COMB}\\
{\htmladdnormallink{A001001}{http://www.research.att.com:80/cgi-bin/access.cgi/as/njas/sequences/eisA.cgi?Anum=001001}}
&
{\htmladdnormallink{A000203}{http://www.research.att.com:80/cgi-bin/access.cgi/as/njas/sequences/eisA.cgi?Anum=000203}}
&-&e\\
{\htmladdnormallink{A001018*}{http://www.research.att.com:80/cgi-bin/access.cgi/as/njas/sequences/eisA.cgi?Anum=001018}}
&
{\htmladdnormallink{A027380*}{http://www.research.att.com:80/cgi-bin/access.cgi/as/njas/sequences/eisA.cgi?Anum=027380}}
&
{full 8-shift}&{POLY}\\
{\htmladdnormallink{A001019*}{http://www.research.att.com:80/cgi-bin/access.cgi/as/njas/sequences/eisA.cgi?Anum=001019}}
&
{\htmladdnormallink{A027381*}{http://www.research.att.com:80/cgi-bin/access.cgi/as/njas/sequences/eisA.cgi?Anum=027381}}
&{full 9-shift}&{POLY}\\
{\htmladdnormallink{A001020*}{http://www.research.att.com:80/cgi-bin/access.cgi/as/njas/sequences/eisA.cgi?Anum=001020}}
&
{\htmladdnormallink{A032166}{http://www.research.att.com:80/cgi-bin/access.cgi/as/njas/sequences/eisA.cgi?Anum=032166}}
&{full 11-shift}&{NECK(11)}\\
{\htmladdnormallink{A001021*}{http://www.research.att.com:80/cgi-bin/access.cgi/as/njas/sequences/eisA.cgi?Anum=001021}}
&
{\htmladdnormallink{A032167}{http://www.research.att.com:80/cgi-bin/access.cgi/as/njas/sequences/eisA.cgi?Anum=032167}}
&{full 12-shift}&{NECK(12)}\\
{\htmladdnormallink{A001022*}{http://www.research.att.com:80/cgi-bin/access.cgi/as/njas/sequences/eisA.cgi?Anum=001022}}
&
{\htmladdnormallink{A060216}{http://www.research.att.com:80/cgi-bin/access.cgi/as/njas/sequences/eisA.cgi?Anum=060216}}
&{full 13-shift}&{NECK(13)}\\
{\htmladdnormallink{A001023*}{http://www.research.att.com:80/cgi-bin/access.cgi/as/njas/sequences/eisA.cgi?Anum=001023}}
&
{\htmladdnormallink{A060217}{http://www.research.att.com:80/cgi-bin/access.cgi/as/njas/sequences/eisA.cgi?Anum=060217}}
&{full 14-shift}&{NECK(14)}\\
{\htmladdnormallink{A001024*}{http://www.research.att.com:80/cgi-bin/access.cgi/as/njas/sequences/eisA.cgi?Anum=001024}}
&
{\htmladdnormallink{A060218}{http://www.research.att.com:80/cgi-bin/access.cgi/as/njas/sequences/eisA.cgi?Anum=060218}}
&{full 15-shift}&{NECK(15)}\\
{\htmladdnormallink{A001025*}{http://www.research.att.com:80/cgi-bin/access.cgi/as/njas/sequences/eisA.cgi?Anum=001025}}
&
{\htmladdnormallink{A060219}{http://www.research.att.com:80/cgi-bin/access.cgi/as/njas/sequences/eisA.cgi?Anum=060219}}
&{full 16-shift}&{NECK(16)}\\
{\htmladdnormallink{A001026*}{http://www.research.att.com:80/cgi-bin/access.cgi/as/njas/sequences/eisA.cgi?Anum=001026}}
&
{\htmladdnormallink{A060220}{http://www.research.att.com:80/cgi-bin/access.cgi/as/njas/sequences/eisA.cgi?Anum=060220}}
&{full 17-shift}&{NECK(17)}\\
{\htmladdnormallink{A001027*}{http://www.research.att.com:80/cgi-bin/access.cgi/as/njas/sequences/eisA.cgi?Anum=001027}}
&
{\htmladdnormallink{A060221}{http://www.research.att.com:80/cgi-bin/access.cgi/as/njas/sequences/eisA.cgi?Anum=060221}}
&
{full 18-shift}&{NECK(18)}\\
{\htmladdnormallink{A001029*}{http://www.research.att.com:80/cgi-bin/access.cgi/as/njas/sequences/eisA.cgi?Anum=001029}}
&
{\htmladdnormallink{A060222}{http://www.research.att.com:80/cgi-bin/access.cgi/as/njas/sequences/eisA.cgi?Anum=060222}}
&
{full 19-shift}&{NECK(19)}\\
{\htmladdnormallink{A001157}{http://www.research.att.com:80/cgi-bin/access.cgi/as/njas/sequences/eisA.cgi?Anum=001157}}
&
{\htmladdnormallink{A000027}{http://www.research.att.com:80/cgi-bin/access.cgi/as/njas/sequences/eisA.cgi?Anum=000027}}
&{-}&{e}\\
{\htmladdnormallink{A001158}{http://www.research.att.com:80/cgi-bin/access.cgi/as/njas/sequences/eisA.cgi?Anum=001158}}
&
{\htmladdnormallink{A000290*}{http://www.research.att.com:80/cgi-bin/access.cgi/as/njas/sequences/eisA.cgi?Anum=000290}}
&{-}&{e}\\
{\htmladdnormallink{A001641?}{http://www.research.att.com:80/cgi-bin/access.cgi/as/njas/sequences/eisA.cgi?Anum=001641}}
&
{\htmladdnormallink{A060166}{http://www.research.att.com:80/cgi-bin/access.cgi/as/njas/sequences/eisA.cgi?Anum=060166}}
&{-}&{-}\\
{\htmladdnormallink{A001642?}{http://www.research.att.com:80/cgi-bin/access.cgi/as/njas/sequences/eisA.cgi?Anum=001642}}
&
{\htmladdnormallink{A060167}{http://www.research.att.com:80/cgi-bin/access.cgi/as/njas/sequences/eisA.cgi?Anum=060167}}
&
{-}&{-}\\
{\htmladdnormallink{A001643?}{http://www.research.att.com:80/cgi-bin/access.cgi/as/njas/sequences/eisA.cgi?Anum=001643}}
&
{\htmladdnormallink{A060168}{http://www.research.att.com:80/cgi-bin/access.cgi/as/njas/sequences/eisA.cgi?Anum=060168}}
&{-}&{-}\\
{\htmladdnormallink{A001700}{http://www.research.att.com:80/cgi-bin/access.cgi/as/njas/sequences/eisA.cgi?Anum=001700}}
&
{\htmladdnormallink{A022553}{http://www.research.att.com:80/cgi-bin/access.cgi/as/njas/sequences/eisA.cgi?Anum=022553}}
&{-}&{-}\\
{\htmladdnormallink{A001945}{http://www.research.att.com:80/cgi-bin/access.cgi/as/njas/sequences/eisA.cgi?Anum=001945}}
&
{\htmladdnormallink{A060169}{http://www.research.att.com:80/cgi-bin/access.cgi/as/njas/sequences/eisA.cgi?Anum=060169}}
&{auto of $\mathbb T^3$}&{-}\\
{\htmladdnormallink{A004146*}{http://www.research.att.com:80/cgi-bin/access.cgi/as/njas/sequences/eisA.cgi?Anum=004146}}
&
{\htmladdnormallink{A032170}{http://www.research.att.com:80/cgi-bin/access.cgi/as/njas/sequences/eisA.cgi?Anum=032170}}
&{auto of $\mathbb T^2$}&{CHK}\\
{\htmladdnormallink{A005809*}{http://www.research.att.com:80/cgi-bin/access.cgi/as/njas/sequences/eisA.cgi?Anum=005809}}
&
{\htmladdnormallink{A060170}{http://www.research.att.com:80/cgi-bin/access.cgi/as/njas/sequences/eisA.cgi?Anum=060170}}
&{-}&{COMB}\\
{\htmladdnormallink{A006953?}{http://www.research.att.com:80/cgi-bin/access.cgi/as/njas/sequences/eisA.cgi?Anum=006953}}
&
{\htmladdnormallink{A060171}{http://www.research.att.com:80/cgi-bin/access.cgi/as/njas/sequences/eisA.cgi?Anum=060171}}
&{-}&{KUMMER?}\\
{\htmladdnormallink{A006954?}{http://www.research.att.com:80/cgi-bin/access.cgi/as/njas/sequences/eisA.cgi?Anum=006954}}
&
{\htmladdnormallink{A060479}{http://www.research.att.com:80/cgi-bin/access.cgi/as/njas/sequences/eisA.cgi?Anum=060479}}
&
{-}&{KUMMER?}\\
{\htmladdnormallink{A011557*}{http://www.research.att.com:80/cgi-bin/access.cgi/as/njas/sequences/eisA.cgi?Anum=011557}}
&
{\htmladdnormallink{A032165*}{http://www.research.att.com:80/cgi-bin/access.cgi/as/njas/sequences/eisA.cgi?Anum=032165}}
&{full 10-shift}&{NECK(10)}\\
{\htmladdnormallink{A023890}{http://www.research.att.com:80/cgi-bin/access.cgi/as/njas/sequences/eisA.cgi?Anum=023890}}
&
{\htmladdnormallink{A005171}{http://www.research.att.com:80/cgi-bin/access.cgi/as/njas/sequences/eisA.cgi?Anum=005171}}
&{-}&{e}\\
{\htmladdnormallink{A027306*}{http://www.research.att.com:80/cgi-bin/access.cgi/as/njas/sequences/eisA.cgi?Anum=027306}}
&
{\htmladdnormallink{A060172}{http://www.research.att.com:80/cgi-bin/access.cgi/as/njas/sequences/eisA.cgi?Anum=060172}}
&{-}&{COMB}\\
{\htmladdnormallink{A035316}{http://www.research.att.com:80/cgi-bin/access.cgi/as/njas/sequences/eisA.cgi?Anum=035316}}
&
{\htmladdnormallink{A010052*}{http://www.research.att.com:80/cgi-bin/access.cgi/as/njas/sequences/eisA.cgi?Anum=010052}}
&{-}&{e}\\
{\htmladdnormallink{A047863*}{http://www.research.att.com:80/cgi-bin/access.cgi/as/njas/sequences/eisA.cgi?Anum=047863}}
&
{\htmladdnormallink{A060224}{http://www.research.att.com:80/cgi-bin/access.cgi/as/njas/sequences/eisA.cgi?Anum=060224}}
&{-}&{-}\\
{\htmladdnormallink{A048578}{http://www.research.att.com:80/cgi-bin/access.cgi/as/njas/sequences/eisA.cgi?Anum=048578}}
&
{\htmladdnormallink{A060477}{http://www.research.att.com:80/cgi-bin/access.cgi/as/njas/sequences/eisA.cgi?Anum=060477}}
&
{4-shift $\cup$ singleton}&{-}\\
{\htmladdnormallink{A056045}{http://www.research.att.com:80/cgi-bin/access.cgi/as/njas/sequences/eisA.cgi?Anum=056045}}
&
{\htmladdnormallink{A060173}{http://www.research.att.com:80/cgi-bin/access.cgi/as/njas/sequences/eisA.cgi?Anum=060173}}
&{-}&{COMB}\\
{0,2,0,6,0,8,0,14,...}
&
{\htmladdnormallink{A000035}{http://www.research.att.com:80/cgi-bin/access.cgi/as/njas/sequences/eisA.cgi?Anum=000035}}
&{-}&{e}\\
{\htmladdnormallink{A059928}{http://www.research.att.com:80/cgi-bin/access.cgi/as/njas/sequences/eisA.cgi?Anum=059928}}
&
{\htmladdnormallink{A060478}{http://www.research.att.com:80/cgi-bin/access.cgi/as/njas/sequences/eisA.cgi?Anum=060478}}
&
{auto of $\mathbb T^{10}$}&{e}\\
{\htmladdnormallink{A059990}{http://www.research.att.com:80/cgi-bin/access.cgi/as/njas/sequences/eisA.cgi?Anum=059990}}
&
{\htmladdnormallink{A060480}{http://www.research.att.com:80/cgi-bin/access.cgi/as/njas/sequences/eisA.cgi?Anum=060480}}
&
{S(1)}&{-}\\
{\htmladdnormallink{A059991}{http://www.research.att.com:80/cgi-bin/access.cgi/as/njas/sequences/eisA.cgi?Anum=059991}}
&
{\htmladdnormallink{A060481}{http://www.research.att.com:80/cgi-bin/access.cgi/as/njas/sequences/eisA.cgi?Anum=060481}}
&{S(2)}&{-}\\
\hline
\caption{\label{table}Exactly realizable sequences.}
\end{longtable}

%\bibliographystyle{plain}
%\bibliography{../bib/refs}

\end{document}